# Gain-Induced Oscillations in Blood Pressure

Roselyn M. Abbiw-Jackson*   William F. Langford †

September 14, 1997

**Abstract**

"Mayer waves" are long-period (6 to 12 seconds) oscillations in arterial blood pressure, which have been observed and studied for more than 100 years in the cardiovascular system of humans and other mammals. A mathematical model of the human cardiovascular system is presented, incorporating parameters relevant to the onset of Mayer waves. The model is analyzed using methods of Lyapunov stability and Hopf bifurcation theory. The analysis shows that increase in the gain of the baroreflex feedback loop controlling venous volume may lead to the onset of oscillations, while changes in the other parameters considered do not affect stability of the equilibrium state. The results agree with clinical observations of Mayer waves in human subjects, both in the period of the oscillations and in the observed age-dependence of Mayer waves. This leads to a proposed explanation of their occurrence, namely that Mayer waves are a "gain-induced instability".

# 1  INTRODUCTION

The existence of fluctuations in blood pressure has been know since the introduction of the recording manometer by C. Ludwig, see [18]. These fluctuations, usually referred to as waves, are classified by various methods including

*Department of Mathematics, University of Maryland, College Park MD, USA 20742

†Department of Mathematics and Statistics, University of Guelph, Guelph ON, Canada N1G 2W1, and The Fields Institue for Research in Mathematical Sciences, 222 College Street, Toronto ON, Canada M5T 3J1. Research supported by NSERC, Canada. Email: wlangfor@uoguelph.ca

the name of the discoverer, the origin, physiological cause, or the frequency or period. The term Mayer waves refers to periodic fluctuations in blood pressure which are slower than respiration in animals with normal respiratory movements. They were announced by S. Mayer in 1876 [15], and hence the name. They are also known as third order waves.

The frequencies reported by various authors for Mayer waves differ considerably [18]. Those described by Mayer in rabbits had a frequency of 6-9 waves/min., while other researchers have found waves with frequencies ranging from 7-12 waves/min. in humans [18]. Some researchers have proposed to designate these waves as the "10-second-rhythm" [18]. The onset of Mayer waves may result in serious physiological implications, such as fainting. Mayer waves are of interest to researchers seeking to fully understand the functioning of the cardiovascular system.

It is generally conceded that Mayer waves appear most often when the subjects are exposed to abnormal conditions [2]. Lack of oxygen, the effects of severe haemorrhage, and other extreme or sudden changes in blood supply to parts of the body favour the appearance of these slow periodic fluctuations [2]. Experiments have shown that when the blood pressure is measured for subjects lying in a supine position and then in a tilted position, there may exist Mayer-like oscillations for the tilted position. A remarkable feature observed in these experiments is that the Mayer waves occur more frequently in younger subjects, and disappear with age [6, 10, 12, 14, 16, 17]. The origin of Mayer waves, however, remains an unsolved problem [18]. At a conference on Mayer waves held in Prague in 1977, various theories were proposed to explain the origin of these waves. Four main theories, the myogenic theory, the central theory, the feedback theory, and the resonance theory were given.

The myogenic theory states that third order waves are due to the inherent property of peripheral systems, (vascular smooth muscle), to exhibit spontaneous rhythmic activity. The central theory postulates that cardiovascular centres in the brain stem generates the slow rhythms in a way similar to the respiratory centre. The feedback theory, as its name suggests, attributes the origin of Mayer waves to delays and nonlinearities in the body's feedback control mechanisms. The resonance theory on the other hand, postulates that one or more of the above factors enable the overall circulatory system to resonate at certain frequencies. Among the four, the feedback theory and the central theory are considered most probable [11, 21, 22]. However, the stability of the central oscillator is questionable [21]. Thus among the two most popular theories, the feedback theory remains the most probable explanation for the existence of Mayer waves.

It is essential, for survival, that blood pressure be controlled to stay within a narrow, safe range. This function is performed by the body's control mechanisms, the fastest being the baroreflex. The baroreceptors are stretch sensors located in the systemic arteries which detect changes in blood pressure. The baroreflex feedback loops respond to baroreceptor impulses to control blood pressure via three mechanisms: heart rate, systemic capillary resistance and venous volume. All three mechanisms are explored in this work.

DeBoer et. al.[5] proposed that blood pressure variability is caused by a time delay in the baroreceptor loop. If this were the case, one would expect the delay to increase with age due to a slowing of the body's responses. This would then cause an increase in the existence of Mayer waves in older people. While we do not discount the contribution of delay, it appears not to be the main cause of Mayer waves. We hypothesize that blood pressure variability may be attributed to a change in feedback gain, extending the work of Wesseling et. al.[21][22]. Previous studies of feedback-control systems in physiology (Glass and Mackey model) [7][13], and in engineering (Watt s regulator model) [8, 19], have shown that an increase in feedback gain can cause a system to change behaviour from a steady state to an oscillating state. This may be called a "gain-induced instability" and has been studied by use of the Hopf Bifurcation Theorem. Since young adult humans tend to have quicker reflexes and better muscle tone than the elderly, they can be expected to have higher gain in the baroreceptor loop. Thus our hypothesis that Mayer waves may be a gain-induced instability is consistent with the observed age-related data.

A mathematical model can be used to give greater insights into the roles of the various mechanisms affecting Mayer waves. Thus the primary objective of this study is to develop a dynamical model for the mammalian circulatory system and use it to analyze blood pressure variability dynamics, as a function of the physiological parameters in the model. The dynamical model is a generalization of the steady-state model of Hoppensteadt and Peskin [9]. In addition to incorporating temporal dynamics, this model will allow investigation of the effects of each of the three baroreflex feedback loops, independently of the other two. This type of experiment has been performed on animals, eg. dogs [20], but is difficult to carry out on human subjects. Parameter values in the mathematical model are chosen to correspond to a typical adult human being.

# 2 MODELLING

While it is desirable to include the behaviour of each cardiovascular component in a model of the circulation, certain components can be lumped together without sacrificing the qualitative behaviour of the system [3, 14]. This section presents the modelling assumptions and the development of the model, first for the basic fluid flow of blood in the cardiovascular system, then with the nonlinear baroreflex control.

## 2.1 MODELLING ASSUMPTIONS

The assumptions and simplifications underlying the mathematical model are stated in the following.

1. The cardiovascular system is a closed-looped hydrodynamic system comprising two heart pumps (the left and right sides of the heart), two large arteries, two veins and the two capillary networks, corresponding to the systemic and pulmonary circulations respectively. The total blood volume is constant in time.

2. The large arteries and veins and the heart are compliance vessels [9], that is, volume is proportional to pressure in these vessels. On the other hand, the smaller arteries and veins in the capillary networks are resistance vessels, that is, flow is proportional to pressure. The unstressed volume of blood vessels is negligible at all parts of the circulation except in the systemic veins.

3. Flow from the heart is continuous, that is the pulsatile nature of blood pressure is neglected and only average pressures and volumes, over the period of the pulse, are dealt with.

4. The pressure in the heart when relaxed, is equal to that of the veins supplying blood to it. That is, the pressure in the right and left hearts are those of the systemic and pulmonary veins respectively. On its expansion stroke (diastole) the heart receives a volume of blood proportional to this venous pressure. The heart pumps out at each contraction stroke (systole) the amount of blood received from the veins on the previous diastole. This is the Frank-Starling model of the heart [22].

5. Cardiovascular blood pressure is controlled by a baroreceptor feedback mechanism, which acts on systemic venous unstressed volume, systemic venous compliance, systemic capillary resistance and heart rate, to counteract changes in systemic arterial pressure.

6. The baroreceptor feedback gain and its dependence on systemic arterial pressure is modelled as a Hill function (described below).

7. Changes in venous volume, systemic resistance and heart rate act independently in parallel on blood pressure.

8. Compliance is constant in all parts of the cardiovascular system except the systemic veins, where it may be varied by the baroreflex.

9. Resistance is constant for the capillary networks of the pulmonary circulation, but may be varied by the baroreflex in the systemic circulation.

## 2.2 MODEL DEVELOPMENT

A notational convention adopted throughout this model is that dynamic variables are represented by lower case letters, while parameters and labels are written in upper case. We first construct a simple linear model of the cardiovascular system, then add the baroreflex control system.

### 2.2.1 Linear Cardiovascular Model

The following linear relationship between volume, $v$, pressure, $p$ and compliance, $C$, in the large vessels (arteries and veins) of the circulation [9], is the mathematical form of Assumption 2. (i.e. that these are compliance vessels):

$$v = C * p \qquad (1)$$

However, Equation (1) suggests that if $p = 0$, then $v = 0$, which is not the case. A more realistic relation is:

$$v = V_D + (C * p) \tag{2}$$

where $V_D$ is the unstressed volume, that is the volume of the vessels at $p = 0$. The unstressed volume in the systemic venous circulation is very important, as over one-half of the venous volume is unstressed volume [4]. (see Assumption 2.).

The flow, $q$ in the vessels of the capillary networks, modelled as in [9] as resistance vessels, is (by Assumption 2):

$$q = \frac{p_A - p_V}{R} \tag{3}$$

Here $p_A$, $p_V$ and $R$ represent the pressure in the arteries and veins and resistance respectively, in either the systemic or pulmonary circulation.

From the Frank-Starling Assumption 4, the following relations for the left and right cardiac outputs, $q_L$ and $q_R$, respectively are given:

$$q_L = F * C_L * p_{PV} = K_L * p_{PV} \tag{4}$$

$$q_R = F * C_R * p_{SV} = K_R * p_{SV} \tag{5}$$

$C_L$ and $C_R$ are the compliances in the left and right hearts respectively and $F$ is the heart beat frequency. The subscripts $A$ and $V$ represent arteries and veins respectively, while the subscripts $S$ and $P$ stand for the systemic and pulmonary circulations, respectively.

The rate of change of volume of an incompressible fluid in a vessel is the difference between the rates of flow of the fluid, into and out of the vessel. Hence, the following differential equations are obtained, for the change of volume of blood in the systemic arteries, systemic veins, pulmonary arteries and pulmonary veins respectively:

$$\frac{dv_{SA}}{dt} = q_L - q_S \tag{6}$$

$$\frac{dv_{SV}}{dt} = q_S - q_R \tag{7}$$

$$\frac{dv_{PA}}{dt} = q_R - q_P \tag{8}$$

$$\frac{dv_{PV}}{dt} = q_P - q_L \tag{9}$$

$q_L$, $q_R$, $q_S$ and $q_P$, represent flow through the left heart, the right heart, the systemic capillaries and the pulmonary capillaries, respectively. The algebraic Equations (1)-(5) may be used to eliminate the flow variables $q$ and the pressure variables $p$ from the system of differential Equations (6)-(9). The result is a system of four differential equations in the four volume variables.

From Assumption 1, the following is obtained:

$$v_{SA} + v_{SV} + v_{PA} + v_{PV} = V_O \tag{10}$$

where $V_O$ is the total blood volume, a constant. Equation (10) indicates that the four volume variables are not independent. As $v_{PA}$ has the smallest value, it is chosen for elimination and a system of three equations in $v_{SA}$, $v_{SV}$ and $v_{PV}$, is obtained (mathematically any one of the four volume variables could be eliminated). Hence the mathematical model of the cardiovascular system consists of the following system of three ordinary differential equations:

$$\frac{dv_{SA}}{dt} = -\frac{v_{SA}}{R_S * C_{SA}} + \frac{v_{SV}}{R_S * C_{SV}} + \frac{F * C_L * v_{PV}}{C_{PV}} - \frac{v_D}{R_S * C_{SV}} \tag{11}$$

$$\frac{dv_{SV}}{dt} = \frac{v_{SA}}{R_S * C_{SA}} - \left(\frac{1}{R_S * C_{SV}} + \frac{F * C_R}{C_{SV}}\right)(v_{SV} - v_D) \tag{12}$$

$$\frac{dv_{PV}}{dt} = \frac{V_O - v_{SA} - v_{SV} - v_{PV}}{R_P * C_{PA}} - \left(\frac{1}{R_P * C_{PV}} + \frac{F * C_L}{C_{PV}}\right)v_{PV} \tag{13}$$

Note that this model is linear in the three volume variables. The system becomes nonlinear, on the inclusion of the baroreflex control loop.

### 2.2.2 Model with Baroreceptor Control

The Hill function is defined by:

$$y = f_n(x) = \frac{x^n}{a^n + x^n} \qquad (14)$$

$$f_n(x) : [0, \infty) \to [0, 1)$$

where $a$ is a constant which corresponds to the point of half the maximum output, $f_n(a) = \frac{1}{2}$. As $n$ increases the curve approaches a perfect "switch" or a step function at $x = a$. The baroreceptor response curve described in the literature strongly resembles a Hill function and therefore is modelled in this paper as:

$$B_n(p_{SA}) = \frac{(p_{SA})^n}{(p_C)^n + (p_{SA})^n} \qquad (15)$$

where $B$ is the total baroreceptor afferent activity, $n$ is a measure of the baroreflex gain, and $p_C$, is the critical arterial pressure. The term "gain" normally is used to represent a ratio of the change in output to a change in input, for very small changes. This however is essentially the mathematical definition of a derivative. Thus, for our model, using the Hill function for the baroreflex response, we define gain by the derivative

$$\mu = \frac{dB}{dp_{SA}} \qquad (16)$$

that is, gain $\mu$ is equal to the slope of the response function (for fixed $n$), at a particular point. For simplicity we take the value of this $\mu$ when $p_{SA} = p_C$, which is precisely

$$\mu = \frac{n}{4} \; . \qquad (17)$$

This is a good measure of the gain as it is very close to the maximum value of the slope (from calculus).

Also, note that $p_{SA}$ is proportional to $v_{SA}$ (since $C_{SA}$ is a constant), so $B$ can be expressed in terms of $v_{SA}$ rather than $p_{SA}$. This yields:

$$B_n(p_{SA}) = B_n(\frac{v_{SA}}{C_{SA}}) = \frac{(v_{SA})^{4\mu}}{(V_C)^{4\mu} + (v_{SA})^{4\mu}} \qquad (18)$$

where $V_C$ is the volume at the critical pressure. It should be noted that other functions with a similar graph shape could be used for the baroreceptor response function, for example the hyperbolic tangent function.

Changes in heart rate, $F$ and systemic capillary resistance, $R_S$ must be in an opposite direction to a change in arterial blood pressure, in order to restore normal pressure. Thus a simple model of the baroreflex action on $F$ is:

$$F = F_0(1 - B) + = \frac{F_0(v_C)^{4\mu}}{(V_C)^{4\mu} + (v_{SA})^{4\mu}} \tag{19}$$

where $F_O$ is a constant. Equation (19) implies that if $B$ approaches 1 (i.e. very large pressure $p_{SA}$), then $F$ will be zero. However, one would expect that in reality $F$ will have a non zero minimum value even when $B$ approaches 1. The following is a more realistic representation:

$$F = F_1(1 - B) + F_2 = \frac{F_1(v_C)^{4\mu}}{(V_C)^{4\mu} + (v_{SA})^{4\mu}} + F_2 \tag{20}$$

where $F_1$ and $F_2$ are constants, and $F_2$ is the value of $F$ when $B = 1$. Similarly, the baroreflex action on $R_S$ is modelled as:

$$R_S = R_1(1 - B) + R_2 = \frac{R_1(v_C)^{4\mu}}{(V_C)^{4\mu} + (v_{SA})^{4\mu}} + R_2 \tag{21}$$

On the other hand, changes in systemic venous compliance, $C_{SV}$, and systemic venous unstressed volume, $V_D$, are in the same direction as a change in arterial blood pressure. Thus the action of the baroreflex on each of $V_D$ and $C_{SV}$ is modelled as:

$$V_D = (D_1 * B) + D_2 = \frac{D_1(v_{SA})^{4\mu}}{(V_C)^{4\mu} + (v_{SA})^{4\mu}} + D_2 \tag{22}$$

$$C_{SV} = (C_1 * B) + C_2 = \frac{C_1(v_{SA})^{4\mu}}{(v_C)^{4\mu} + (v_{SA})^{4\mu}} + C_2 \tag{23}$$

Now, from the basic linear cardiovascular model (11)(12)(13), four different nonlinear cardiovascular models are obtained, corresponding to insertion of the baroreceptor feedback function into each of F, $R_S$, $C_{SV}$ and $V_D$, as above. This allows independent investigations of each of the four feedback

Table 1: Typical parameter values for an adult human being(Hoppensteadt and Peskin, 1992).

| PARAMETER | NORMAL VALUE |
|---|---|
| Compliance in Systemic Arteries, $C_{SA}$ | 0.01 litres/mm Hg |
| Compliance in Pulmonary Arteries, $C_{PA}$ | 0.00667 litres/mm Hg |
| Compliance in Pulmonary Veins, $C_{PV}$ | 0.08 litres/mm Hg |
| Systemic resistance, $R_S$ | 17.5 mm Hg/(litre/min.) |
| Pulmonary resistance, $R_P$ | 1.79 mm Hg/(litre/min.) |
| Compliance in Right Heart, $C_R$ | 035 litres/mm Hg |
| Compliance in Left Heart, $C_L$ | 0.014 litres/mm Hg |
| Heart rate, $F$ | 80 beats/min. |

loops, which would be very difficult to carry out in experiments on live subjects.

## 2.3 PARAMETER DETERMINATION

Many of the parameters in this model are available in the literature, as displayed in Table 1. The remaining parameters, $F_1$, $F_2$, $R_1$, $R_2$, $C_1$, $C_2$, $D_1$, $D_2$, $C_{SV}$, $V_C$ and $V_D$ are not found in the literature and need to be determined.

### 2.3.1 Critical Volume, $V_C$

No value of $V_C$ is in the literature. However the normal resting value of $v_{SA}$ is known to be 1.0 litre. It is assumed that the resting and critical states are the same and hence $V_C$ is taken as 1.0 litre.

### 2.3.2 Normal Systemic Venous Unstressed Volume, $V_D$

An exact normal value of $V_D$ is not found in the literature. However, Coleman (1985), gives the value of $V_D$ as "over half" the systemic venous volume. As the normal value of $v_{SV}$ is 3.5 litres, the normal value of $V_D$ is taken as 2.0 litres in this study.

### 2.3.3 Normal Systemic Venous Compliance, $C_{SV}$

The normal systemic venous compliance is given as 1.75 litres/mm Hg (Hoppensteadt and Peskin, 1992). This value of $C_{SV}$ does not account for systemic venous unstressed volume. However, this study considers the systemic venous unstressed volume. Using Equation 2, with $v_{SV} = 3.5$ litres, $V_D = 2.0$ litres and $p_{SV} = 2$ mm Hg, the value of $C_{SV}$ is computed as 0.75 litres/mm Hg.

### 2.3.4 Normalized Hill Function Constants, $F_1, F_2, R_1, R_2, C_1, C_2, D_1,$ and $D_2$

These constants are required for the use of the Hill function to model the baroreceptor afferent activity, $B$, for the systemic venous compliance, systemic venous unstretched volume, systemic resistance and heart rate. This current study appears to be the first time such an approach has been taken, which is why these constants are not available. Therefore, sets of different values of each of these constants are investigated in this study, over ranges which yield baroreceptor responses consistent with experimental observations.

## 3 ANALYSIS AND RESULTS

The control of the baroreflex on heart rate, $F$, systemic capillary resistance, $R_S$, systemic venous unstressed volume, $V_D$, and systemic venous compliance, $C_{SV}$, are investigated individually in the mathematical model. The nonlinear baroreflex response function is substituted in each parameter in turn, to obtain the corresponding model for investigating the baroreflex effect on the parameter under consideration. See the Appendixes and [1] for details. We analyzed each of the four models to find out if a bifurcation occurs as the baroreceptor gain $\mu$ varies, using Hopf's Bifurcation Theorem. First, the steady state solution of each model was found. Then the system of equations was linearised at this steady state. The value of the Jacobian matrix obtained at the steady state was found. Since it is a real $3 \times 3$ matrix, with constant real entries, the eigenvalues are either all three real or else one real and two complex conjugate. The eigenvalues of the resulting matrix were found. In each case there existed, for some values of $\mu$, one real and two complex eigenvalues. The real part of the pair of complex conjugate eigenvalues was plotted as a function of $\mu$, to find out if a crossing point existed. The value of $\mu$ at which the real part crosses the $\mu$-axis is what is known as the crossing point or Hopf bifurcation point. At this point there exists a pair of purely imaginary eigenvalues $\pm i\omega$ and the steady state is said to be nonhyperbolic. When a crossing point was found, the imaginary part $\omega$ of the complex eigenvalues was plotted to obtain its value at the crossing point. The third (real) eigenvalue always remained negative. According to the general theory of Liapunov stability, when the real part of the complex eigenvalues crosses from negative to positive, the equilibrium state changes from asymptotically stable to unstable. The computation of the eigenvalues and the plotting of the curves was done using Maple. From the Hopf Bifurcation Theorem [8], generically at such a crossing point, a periodic solution is either created or destroyed. Further numerical computations verify the existence of a stable limit cycle near the crossing point. The imaginary part $\omega$ at the crossing point gives a good approximation to the frequency of the resulting oscillations.

## 3.1 BAROREFLEX CONTROL OF HEART RATE

Models with $R_{\text{S}}$, $C_{\text{SV}}$ and $V_{\text{D}}$ taken as constants and $F$ given by Equation 20 were considered. Assuming a normal heart rate of 80 beats/min., values of $F_1$ and $F_2$ considered were: $F_1 = 160$ beats/min. and $F_2 = 0$ beats/min., $F_1 = 80$ beats/min. and $F_2 = 40$ beats/min., and $F_1 = 40$ beats/min. and $F_2 = 60$ beats/min. All of these models exhibited a stable steady-state, for all values of gain $\mu$ tested. No evidence of waves was found. Details of the models are given in [1].

## 3.2 BAROREFLEX CONTROL OF SYSTEMIC RESISTANCE

Models with the baroreflex affecting only systemic resistance $R_{\text{S}}$, while $C_{\text{SV}}$, $V_{\text{D}}$ and $F$ are taken as constants, are considered next. $R_{\text{S}}$ is given by Equation 21 and a typical value of systemic resistance is 17.5 mmHg/(litre/min.). Values of $R_1$ and $R_2$ used are: $R_1 = 35$ mm Hg/(litre/min.) and $R_2 = 0$ mm Hg/(litre/min.), $R_1 = 20$ mm Hg/(litre/min.) and $R_2 = 7.5$ mm Hg/(litre/min.), and $R_1 = 15$ mm Hg/(litre/min.) and $R_2 = 10$ mm Hg/(litre/min.). These models had a stable steady-state for all values of gain, $\mu$ tested, and showed no indications of waves. Details of the models are given in [1].

## 3.3 BAROREFLEX CONTROL OF VENOUS VOLUME

The baroreflex influences the systemic venous volume through the unstressed volume, $V_D$ and the compliance, $C_{SV}$. Models with the baroreflex controlling $C_{SV}$ and $V_D$ individually were considered (see Appendix for equations). For both cases the models exhibited unstable steady-states for gains past a crossing point with pure imaginary eigenvalues. Figure 1 displays four graphs obtained for models with the baroreflex controlling unstressed venous volume. It shows graphs of the real parts of the complex eigenvalues, for models with $D_2$ equal to 0, 0.5, 1.0, and 1.5 litres, respectively. Note that $Re(\lambda)$ crosses through zero in all cases. This implies a Hopf bifurcation, giving birth to an oscillation or wave.

Similarly, Figure 2 is obtained from models with the baroreflex controlling venous compliance only. It shows graphs of the real parts of the complex eigenvalues for three cases of models with $C_2$ equal to 0, 0.25, and 0.5 litres/mm Hg, respectively. All three cases give a Hopf bifurcation.

The values of the imaginary parts of the complex eigenvalues at the crossing points give the angular frequency of the oscillations produced. From these frequencies, the periods of oscillation of all these models were found to be between 7 and 12 secs. Note that this is in perfect agreement with reported values of Mayer waves in human subjects.

Figure 3 shows phase portraits obtained for the system with $D_2 = 0$, for gains $\mu$ equal to 10 and 20. At $\mu = 10$, there exists no limit cycle. However, for $\mu = 20$, a limit cycle exists in confirmation of the above analysis. In all cases considered, when a crossing point existed, a limit cycle was found to be stable and supercritical (that is, for $\mu$ above the crossing point).

# 4 DISCUSSION

Models with only heart rate $F$, or systemic capillary resistance $R_S$, controlled by the baroreflex did not exhibit a Hopf Bifurcation, while the models with systemic venous compliance, $C_{SV}$ or systemic venous unstressed volume, $V_D$ controlled by the baroreflex were capable of Hopf bifurcation. Hence the effect of the baroreflex on $F$ and $R_S$, individually, is not the cause of oscillations. However, if the effect of the baroreflex on $F$ and $R_S$ were combined in models with the effects on $V_D$ or $C_{SV}$, they may play a part in causing instability.

It is observed that for all models with $V_D$ and $C_{SV}$ individually controlled by the baroreflex, the real part of the complex eigenvalues increases as the gain $\mu$ increases and the graph crosses the $\mu$-axis at a positive value of $\mu$. This implies a Hopf bifurcation and the presence of a limit cycle oscillation. The stability of this limit cycle oscillation has been verified numerically. The similarity of the results obtained for models with $C_{SV}$ and $V_D$ individually controlled by the baroreflex is to be expected, as the two have similar effects on blood flow. The model with the baroreflex controlling only $C_{SV}$ is more stable than that with the baroreflex controlling only $V_D$. Thus it appears that the baroreflex control of $V_D$ is more important than the control of $C_{SV}$ where gain induced instability is concerned, but control of $C_{SV}$ also plays an important part. Remarkably, the periods of the oscillations fall within the range of 7 to 12 secs. given by Penaz [18].

Measurements of Wesseling et. al. [21] on human subjects revealed a periodic fluctuation in heart rate which was synchronized with the Mayer waves in blood pressure. In these experiments, the same period was observed for the modulation of the heart rate as for the blood pressure waves. The present model provides a possible explanation for this heart rate modulation. Assume that Mayer waves have arisen through the mechanism of this paper. Then the baroreceptors would sense a fluctuating blood pressure in the systemic arteries, and would therefore exert a fluctuating feedback control, affecting the heart rate in synchronization with the Mayer waves, as observed in [21].

Further insights were obtained on varying two parameters in the model simultaneously; namely, the gain parameter $\mu$ together with either one of $D_2$ or $C_2$ (the minimum unstressed systemic venous volume or compliance respectively). As either of $D_2$ or $C_2$, increases, the value of gain at which the graph crosses the $\mu$-axis increases. This suggests that an increase in $D_2$ or $C_2$ increases the stability of blood pressure. Changing the value of $D_2$ or $C_2$ causes resetting of the baroreceptor curve. This happens in an individual with a time constant of 10 hours and so is not of importance to this investigation. However as different people may have different $D_2$ and $C_2$ values, different people can be expected to take different times before Mayer waves are observed, when subjected to identical Mayer wave inducing stresses. In particular it could be expected that young and old people will have different $D_2$ and $C_2$ values and this may explain the difference in the incidence in Mayer waves observed in young and old people. As we would expect larger $D_2$ and $C_2$ values in older people, corresponding to veins which have become stretched and less fit, our observation that stability increases with increasing values of $D_2$ and $C_2$ is in agreement with the reactions observed experimentally [6, 9, 11, 14, 16]. In Figure 4, $D_2$ is plotted against the crossing point value of gain $\mu$. The top left region represents the parameter values for which the equilibrium state is stable, and corresponds to older subjects, who would tend to have smaller gains $\mu$, and larger $D_2$ values. The lower right region represents unstable equilibria, susceptible to oscillations, and corresponds to youth. Thus, stability depends on both the baroreflex gain $\mu$ and $D_2$. Thus, young adults, with high gain and small $D_2$, are in the unstable region, while older adults with the opposite characteristics are in the stable region. A similar situation holds for $\mu$ and $C_2$. We conclude that the existence of Mayer waves and their disappearance with age can be explained, at least in part, as a case of gain induced instability.

# 5 APPENDIX A: CONTROL OF VENOUS VOLUME

A mathematical model is constructed in which the systemic venous unstressed volume $V_D$ is controlled by the baroreflex, while $R_S$, $C_{SV}$ and $F$ are assumed to remain constant. $V_D$ is given by Equation 22. Different choices of the constants $D_1$ and $D_2$ are considered. The stability of the equilibrium state is investigated, by computation of the eigenvalues of the Jacobian matrix, using Maple.

## 5.1 $D_1 = 4.0$ litres and $D_2 = 0$ litres

The circulation is then described by the following system of equations:

$$\frac{dv_{SA}}{dt} = -\frac{40v_{SA}}{7} + \frac{8v_{SV}}{105} + (14v_{PV}) - \frac{32(v_{SA})^{4\mu}}{105(1+(v_{SA})^{4\mu})}$$

$$\frac{dv_{SV}}{dt} = \frac{40v_{SA}}{7} - \frac{80v_{SV}}{21} + \frac{320(v_{SA})^{4\mu}}{21(1+(v_{SA})^{4\mu})}$$

$$\frac{dv_{PV}}{dt} = 420 - (84v_{SA}) - (84v_{SV}) - (105v_{PV})$$

The steady state values of the system are found to be: $v_{\text{SA}} = 1.0$ litres, $v_{\text{SV}} = 3.5$ litres, and $v_{\text{PV}} = 0.4$ litres. Linearization of the model at the steady state gives the matrix A.

$$A = \begin{pmatrix} \frac{-40}{7} - \frac{32\mu}{105} & \frac{8}{105} & 14 \\ \frac{40}{7} + \frac{320\mu}{21} & -\frac{80}{21} & 0 \\ -84 & -84 & -105 \end{pmatrix}$$

The eigenvalues of A are:

$$\lambda_1 = U^{1/3} - V - \frac{2405}{63} - \frac{32\mu}{315}$$

$$\lambda_2 = \frac{U^{1/3}}{2} + \frac{V}{2} - \frac{2405}{63} - \frac{32\mu}{315} + \frac{i}{2}\sqrt{3}(U^{1/3} + V)$$

$$\lambda_3 = \frac{U^{1/3}}{2} + \frac{V}{2} - \frac{2405}{63} - \frac{32\mu}{315} - \frac{i}{2}\sqrt{3}(U^{1/3} + V)$$

where

$$U = \frac{-5103633725}{250047} - \frac{3618065648\mu}{416745} + \frac{184832\mu^2}{416745} - \frac{32768\mu^3}{31255875} + \frac{8}{2835}(4488999170550 + 45082290017400\mu + 9458754238932\mu^2 - 942218496\mu^3 + 2248704\mu^4)^{1/2}$$

and

$$V = \frac{\frac{-2876953}{3969} + \frac{11552\mu}{3969} - \frac{1024\mu^2}{99225}}{U^{1/3}}$$

For values of $\mu$ of interest, eigenvalue $\lambda_1$ is real and negative. Eigenvalues $\lambda_2$ and $\lambda_3$ are complex conjugates, with real part which crosses through zero from negative to positive as $\mu$ increases, near $\mu = 18$, as shown in Figure 1a. Phase portraits, on each side of the crossing point, are shown in Figure 3.

## 5.2  $D_1 = 3.0$ litres and $D_2 = 0.5$ litres

The circulation is then described by the following system of equations:

$$\frac{dv_{SA}}{dt} = -\frac{40v_{SA}}{7} + \frac{8v_{SV}}{105} + (14v_{PV}) - \frac{28(v_{SA})^{4\mu} + 4}{105(1 + (v_{SA})^{4\mu})}$$

$$\frac{dv_{SV}}{dt} = \frac{40v_{SA}}{7} - \frac{80v_{SV}}{21} + \frac{40(1 + 7(v_{SA})^{4\mu})}{21(1 + (v_{SA})^{4\mu})}$$

$$\frac{dv_{PV}}{dt} = 420 - (84v_{SA}) - (84v_{SV}) - (105v_{PV})$$

The steady state values of the system are found to be: $v_{SA} = 1.0$ litres, $v_{SV} = 3.5$ litres, and $v_{PV} = 0.4$ litres. Linearization of the model at the steady state gives the matrix A.

$$A = \begin{pmatrix} \frac{-40}{7} - \frac{8\mu}{35} & \frac{8}{105} & 14 \\ \frac{40}{7} + \frac{80\mu}{7} & -\frac{80}{21} & 0 \\ -84 & -84 & -105 \end{pmatrix}$$

The eigenvalues of A are:

$$\lambda_1 = U^{1/3} - V - \frac{2405}{63} - \frac{8\mu}{105}$$

$$\lambda_2 = \frac{U^{1/3}}{2} + \frac{V}{2} - \frac{2405}{63} - \frac{8\mu}{105} + \frac{i}{2}\sqrt{3}(U^{1/3} + V)$$

$$\lambda_3 = \frac{U^{1/3}}{2} + \frac{V}{2} - \frac{2405}{63} - \frac{8\mu}{105} - \frac{i}{2}\sqrt{3}(U^{1/3} + V)$$

where

$$U = \frac{-5103633725}{250047} - \frac{904516412\mu}{138915} + \frac{11552\mu^2}{46305} - \frac{512\mu^3}{1157625} + \frac{4}{2835}(17955996682200 + 135246870052200\mu + 21282197037597\mu^2 - 1589993712\mu^3 + 2846016\mu^4)^{1/2}$$

and

$$V = \frac{\frac{-2876953}{3969} + \frac{2888\mu}{1323} - \frac{64\mu^2}{11025}}{U^{1/3}}$$

For values of $\mu$ of interest, eigenvalue $\lambda_1$ is real and negative. Eigenvalues $\lambda_2$ and $\lambda_3$ are complex conjugates, with real part which crosses through zero from negative to positive as $\mu$ increases, near $\mu = 24$, as shown in Figure 1b.

## 5.3 $D_1 = 2.0$ litres and $D_2 = 1.0$ litres

The circulation is then described by the following system of equations:

$$\frac{dv_{SA}}{dt} = -\frac{40v_{SA}}{7} + \frac{8v_{SV}}{105} + (14v_{PV}) - \frac{24(v_{SA})^{4\mu} + 8}{105(1 + (v_{SA})^{4\mu})}$$

$$\frac{dv_{SV}}{dt} = \frac{40v_{SA}}{7} - \frac{80v_{SV}}{21} + \frac{80(1 + 3(v_{SA})^{4\mu})}{21(1 + (v_{SA})^{4\mu})}$$

$$\frac{dv_{PV}}{dt} = 420 - (84v_{SA}) - (84v_{SV}) - (105v_{PV})$$

The steady state values of the system are found to be: $v_{\text{SA}} = 1.0$ litres, $v_{\text{SV}} = 3.5$ litres, and $v_{\text{PV}} = 0.4$ litres. Linearization of the model at the steady state gives the matrix A.

$$A = \begin{pmatrix} \frac{-40}{7} - \frac{16\mu}{105} & \frac{8}{105} & 14 \\ \frac{40}{7} + \frac{160\mu}{21} & -\frac{80}{21} & 0 \\ -84 & -84 & -105 \end{pmatrix}$$

The eigenvalues of A are:

$$\lambda_1 = U^{1/3} - V - \frac{2405}{63} - \frac{16\mu}{315}$$

$$\lambda_2 = \frac{U^{1/3}}{2} + \frac{V}{2} - \frac{2405}{63} - \frac{16\mu}{315} + \frac{i}{2}\sqrt{3}(U^{1/3} + V)$$

$$\lambda_3 = \frac{U^{1/3}}{2} + \frac{V}{2} - \frac{2405}{63} - \frac{16\mu}{315} - \frac{i}{2}\sqrt{3}(U^{1/3} + V)$$

where

$$\begin{aligned} U &= \frac{-5103633725}{250047} - \frac{1809032824\mu}{416745} + \frac{46208\mu^2}{416745} - \frac{4096\mu^3}{31255875} + \\ & \quad \frac{8}{2835}(4488999170550 + 22541145008700\mu + \\ & \quad 2364688559733\mu^2 - 117777312\mu^3 + 140544\mu^4)^{1/2} \end{aligned}$$

and

$$V = \frac{\frac{-2876953}{3969} + \frac{5776\mu}{3969} - \frac{256\mu^2}{99225}}{U^{1/3}}$$

For values of $\mu$ of interest, eigenvalue $\lambda_1$ is real and negative. Eigenvalues $\lambda_2$ and $\lambda_3$ are complex conjugates, with real part which crosses through zero from negative to positive as $\mu$ increases, near $\mu = 36$, as shown in Figure 1c.

## 5.4   $D_1 = 1.0$ litres and $D_2 = 1.5$ litres

The circulation is then described by the following system of equations:

$$\frac{dv_{SA}}{dt} = -\frac{40v_{SA}}{7} + \frac{8v_{SV}}{105} + (14v_{PV}) - \frac{20(v_{SA})^{4\mu} + 12}{105(1+(v_{SA})^{4\mu})}$$

$$\frac{dv_{SV}}{dt} = \frac{40v_{SA}}{7} - \frac{80v_{SV}}{21} + \frac{40(3+5(v_{SA})^{4\mu})}{21(1+(v_{SA})^{4\mu})}$$

$$\frac{dv_{PV}}{dt} = 420 - (84v_{SA}) - (84v_{SV}) - (105v_{PV})$$

The steady state values of the system are found to be: $v_{SA} = 1.0$ litres, $v_{SV} = 3.5$ litres, and $v_{PV} = 0.4$ litres. Linearization of the model at the steady state gives the matrix A.

$$A = \begin{pmatrix} \frac{-40}{7} - \frac{8\mu}{105} & \frac{8}{105} & 14 \\ \frac{40}{7} + \frac{80\mu}{21} & -\frac{80}{21} & 0 \\ -84 & -84 & -105 \end{pmatrix}$$

The eigenvalues of A are:

$$\lambda_1 = U^{1/3} - V - \frac{2405}{63} - \frac{8\mu}{315}$$

$$\lambda_2 = \frac{U^{1/3}}{2} + \frac{V}{2} - \frac{2405}{63} - \frac{8\mu}{315} + \frac{i}{2}\sqrt{3}(U^{1/3} + V)$$

$$\lambda_3 = \frac{U^{1/3}}{2} + \frac{V}{2} - \frac{2405}{63} - \frac{8\mu}{315} - \frac{i}{2}\sqrt{3}(U^{1/3} + V)$$

where

$$U = \frac{-5103633725}{250047} - \frac{904516412\mu}{416745} + \frac{11552\mu^2}{416745} - \frac{512\mu^3}{31255875} +$$
$$\frac{4}{2835}(17955996682200 + 45082290017400\mu +$$
$$2364688559733\mu^2 - 58888656\mu^3 + 35136\mu^4)^{1/2}$$

and

$$V = \frac{\frac{-2876953}{3969} + \frac{2888\mu}{3969} - \frac{64\mu^2}{99225}}{U^{1/3}}$$

For values of $\mu$ of interest, eigenvalue $\lambda_1$ is real and negative. Eigenvalues $\lambda_2$ and $\lambda_3$ are complex conjugates, with real part which crosses through zero from negative to positive as $\mu$ increases, near $\mu = 71$, as shown in Figure 1d.

# 6 APPENDIX B: CONTROL OF VENOUS COMPLIANCE

A mathematical model is constructed in which the systemic venous compliance $C_{SV}$ is controlled by the baroreflex, while $R_s$, $V_D$ and $F$, are held constant. $C_{SV}$ is given by Equation 23. The stability of the model, for different values of $C_1$ and $C_2$ is explored. The results of these calculations are presented in Figure 2.

## 6.1 $C_1 = 1.5$ litres/mm Hg and $C_2 = 0$ litres/mm Hg

The circulation is then described by the following system of equations:

$$\frac{dv_{SA}}{dt} = -\frac{40v_{SA}}{7} + \frac{4v_{SV}(1 + (v_{SA})^{4\mu})}{105(v_{SA})^{4\mu}} + (14v_{PV}) - \frac{8(1 + (v_{SA})^{4\mu})}{105(v_{SA})^{4\mu}}$$

$$\frac{dv_{SV}}{dt} = \frac{40v_{SA}}{7} - \frac{40v_{SV}(1 + (v_{SA})^{4\mu})}{21(v_{SA})^{4\mu}} + \frac{80(1 + (v_{SA})^{4\mu})}{21(v_{SA})^{4\mu}}$$

$$\frac{dv_{PV}}{dt} = 420 - (84v_{SA}) - (84v_{SV}) - (105v_{PV})$$

The steady state values of the system are found to be: $v_{\text{SA}} = 1.0$ litres, $v_{\text{SV}} = 3.5$ litres, and $v_{\text{P}_{\text{V}}} = 0.4$ litres. Linearization of the model at the steady state gives the matrix A.

$$A = \begin{pmatrix} \frac{-40}{7} - \frac{8\mu}{35} & \frac{8}{105} & 14 \\ \frac{40}{7} + \frac{80\mu}{21} & -\frac{80}{21} & 0 \\ -84 & -84 & -105 \end{pmatrix}$$

The eigenvalues of A are:

$$\lambda_1 = U^{1/3} - V - \frac{2405}{63} - \frac{8\mu}{105}$$

$$\lambda_2 = \frac{U^{1/3}}{2} + \frac{V}{2} - \frac{2405}{63} - \frac{8\mu}{105} + \frac{i}{2}\sqrt{3}(U^{1/3} + V)$$

$$\lambda_3 = \frac{U^{1/3}}{2} + \frac{V}{2} - \frac{2405}{63} - \frac{8\mu}{105} - \frac{i}{2}\sqrt{3}(U^{1/3} + V)$$

where

$$\begin{aligned} U &= \frac{-5103633725}{250047} - \frac{904516412\mu}{138915} + \frac{11552\mu^2}{46305} - \frac{512\mu^3}{1157625} + \\ &\quad \frac{4}{2835}(17955996682200 + 135246870052200\mu + \\ &\quad 21282197037597\mu^2 - 1589993712\mu^3 + 2846016\mu^4)^{1/2} \end{aligned}$$

and

$$V = \frac{\frac{-2876953}{3969} + \frac{2888\mu}{1323} - \frac{64\mu^2}{11025}}{U^{1/3}}$$

For values of $\mu$ of interest, eigenvalue $\lambda_1$ is real and negative. Eigenvalues $\lambda_2$ and $\lambda_3$ are complex conjugates, with real part which crosses through zero from negative to positive as $\mu$ increases, near $\mu = 24$, as shown in Figure 2a.

## 6.2 $C_1 = 1.0$ litres/mm Hg and $C_2 = 0.25$ litres/mm Hg

The circulation is then described by the following system of equations:

$$\frac{dv_{SA}}{dt} = -\frac{40v_{SA}}{7} + \frac{8v_{SV}(1+(v_{SA})^{4\mu})}{35(1+5(v_{SA})^{4\mu})} + (14v_{PV}) - \frac{16(1+(v_{SA})^{4\mu})}{35(1+5(v_{SA})^{4\mu})}$$

$$\frac{dv_{SV}}{dt} = \frac{40v_{SA}}{7} - \frac{80v_{SV}(1+(v_{SA})^{4\mu})}{7(1+5(v_{SA})^{4\mu})} + \frac{160(1+(v_{SA})^{4\mu})}{7(1+5(v_{SA})^{4\mu})}$$

$$\frac{dv_{PV}}{dt} = 420 - (84v_{SA}) - (84v_{SV}) - (105v_{PV})$$

The steady state values of the system are found to be: $v_{SA} = 1.0$ litres, $v_{SV} = 3.5$ litres, and $v_{PV} = 0.4$ litres. Linearization of the model at the steady state gives the matrix A.

$$A = \begin{pmatrix} \frac{-40}{7} - \frac{16\mu}{105} & \frac{8}{105} & 14 \\ \frac{40}{7} + \frac{160\mu}{21} & -\frac{80}{21} & 0 \\ -84 & -84 & -105 \end{pmatrix}$$

The eigenvalues of A are:

$$\lambda_1 = U^{1/3} - V - \frac{2405}{63} - \frac{16\mu}{315}$$

$$\lambda_2 = \frac{U^{1/3}}{2} + \frac{V}{2} - \frac{2405}{63} - \frac{16\mu}{315} + \frac{i}{2}\sqrt{3}(U^{1/3} + V)$$

$$\lambda_3 = \frac{U^{1/3}}{2} + \frac{V}{2} - \frac{2405}{63} - \frac{16\mu}{315} - \frac{i}{2}\sqrt{3}(U^{1/3} + V)$$

where

$$U = \frac{-5103633725}{250047} - \frac{1809032824\mu}{416745} + \frac{46208\mu^2}{416745} - \frac{4096\mu^3}{31255875} +$$
$$\frac{8}{2835}(4488999170550 + 22541145008700\mu +$$
$$2364688559733\mu^2 - 117777312\mu^3 + 140544\mu^4)^{1/2}$$

and

$$V = \frac{\frac{-2876953}{3969} + \frac{5776\mu}{3969} - \frac{256\mu^2}{99225}}{U^{1/3}}$$

For values of $\mu$ of interest, eigenvalue $\lambda_1$ is real and negative. Eigenvalues $\lambda_2$ and $\lambda_3$ are complex conjugates, with real part which crosses through zero from negative to positive as $\mu$ increases, near $\mu = 36$, as shown in Figure 2b.

## 6.3   $C_1 = 0.5$ litres/mm Hg and $C_2 = 0.5$ litres/mm Hg and

The circulation is then described by the following system of equations:

$$\frac{dv_{SA}}{dt} = -\frac{40v_{SA}}{7} + \frac{4v_{SV}(1 + (v_{SA})^{4\mu})}{35(1 + 2(v_{SA})^{4\mu})} + (14v_{PV}) - \frac{8(1 + (v_{SA})^{4\mu})}{35(1 + 2(v_{SA})^{4\mu})}$$

$$\frac{dv_{SV}}{dt} = \frac{40v_{SA}}{7} - \frac{40v_{SV}(1 + (v_{SA})^{4\mu})}{7(1 + 2(v_{SA})^{4\mu})} + \frac{80(1 + v_{SA}^{4\mu})}{7(1 + 2(v_{SA})^{4\mu})}$$

$$\frac{dv_{PV}}{dt} = 420 - (84v_{SA}) - (84v_{SV}) - (105v_{PV})$$

The steady state values of the system are found to be: $v_{\text{SA}} = 1.0$ litres, $v_{\text{SV}} = 3.5$ litres, and $v_{\text{PV}} = 0.4$ litres. Linearization of the model at the steady state gives the matrix A.

$$A = \begin{pmatrix} \frac{-40}{7} - \frac{8\mu}{105} & \frac{8}{105} & 14 \\ \frac{40}{7} + \frac{80\mu}{21} & -\frac{80}{21} & 0 \\ -84 & -84 & -105 \end{pmatrix}$$

The eigenvalues of A are:

$$\lambda_1 = U^{1/3} - V - \frac{2405}{63} - \frac{8\mu}{315}$$

$$\lambda_2 = \frac{U^{1/3}}{2} + \frac{V}{2} - \frac{2405}{63} - \frac{8\mu}{315} + \frac{i}{2}\sqrt{3}(U^{1/3} + V)$$

$$\lambda_3 = \frac{U^{1/3}}{2} + \frac{V}{2} - \frac{2405}{63} - \frac{8\mu}{315} - \frac{i}{2}\sqrt{3}(U^{1/3} + V)$$

where

$$\begin{aligned} U &= \frac{-5103633725}{250047} - \frac{904516412\mu}{416745} + \frac{11552\mu^2}{416745} - \frac{512\mu^3}{31255875} + \\ &\quad \frac{4}{2835}(17955996682200 + 45082290017400\mu + \\ &\quad 2364688559733\mu^2 - 58888656\mu^3 + 35136\mu^4)^{1/2} \end{aligned}$$

and

$$V = \frac{\frac{-2876953}{3969} + \frac{2888\mu}{3969} - \frac{64\mu^2}{99225}}{U^{1/3}}$$

For values of $\mu$ of interest, eigenvalue $\lambda_1$ is real and negative. Eigenvalues $\lambda_2$ and $\lambda_3$ are complex conjugates, with real part which crosses through zero from negative to positive as $\mu$ increases, near $\mu = 71$, as shown in Figure 2c.

# LIST OF FIGURES